# Eichler cohomology in general weights using spectral theory

Michael O. Neururer[1]

*Dedicated to the memory of Marvin Knopp.*

**Abstract** In this paper, we construct a pairing between modular forms of positive real weight and elements of certain Eichler cohomology groups that were introduced by Knopp in 1974. We use spectral theory of automorphic forms to show that this pairing is perfect for all positive weights except 1. The approach in this paper gives a new proof of a theorem by Knopp and Mawi from 2010 for all real weights excluding 1 and also a version of this theorem for vector-valued modular forms.

## 1 Introduction

Let $\Gamma \subseteq \mathrm{SL}_2(\mathbb{R})$ be a finitely generated Fuchsian group of the first kind that contains translations and $-I$, where $I$ is the identity matrix. The interpretation of modular forms for $\Gamma$ as elements in certain cohomology groups was first discovered by Eichler [5]. The following theorem is due to him in the case of even weights and trivial multiplier system. The general case was proved later by Gunning [10].

**Theorem 1** *Let $r$ be a non-positive integer, $v$ a weight $2-r$ multiplier system for $\Gamma$ and $P_r$ the vector space of polynomials with coefficients in $\mathbb{C}$ of degree $\leq -r$. Then*

$$M_{2-r}(\Gamma, v) \oplus S_{2-r}(\Gamma, \overline{v}) \cong H^1_{r,v}(P_r).$$

✉ Michael O. Neururer
  neururer@mathematik.tu-darmstadt.de

[1] School of Mathematics, University of Nottingham, University Park, Nottingham NG7 2RD, UK

Here $P_r$ is viewed as a $\Gamma$-module with the $|_{r,v}$ action and $H^1_{r,v}$ is the first cohomology group. Theorem 1 has many applications in the theory of modular forms and the study of critical values of their $L$-functions, e.g. in algebraicity results like Manin's period theorem [16].

The subject of this paper is a variant of Theorem 1 in the case of arbitrary real weight. Knopp first formulated it in 1974 [12]. To a cusp form $g$ of real weight $2 - r$ and multiplier system $\overline{v}$, he associated a cocycle with values in a space of functions $\mathcal{P}$ by

$$\phi_g(z) : \gamma \mapsto \left[ \int_{\gamma^{-1}\infty}^{\infty} g(\tau)(\tau - \overline{z})^{-r} d\tau \right]^{-}.$$

This induces an injective map from $S_{2-r}(\Gamma, \overline{v})$ to $H^1_{r,v}(\Gamma, \mathcal{P})$ (see Corollary 2) and he conjectured that this map is actually an isomorphism, but was only able to prove this for the cases $r \leq 0$ and $r \geq 2$. In the case $r < 0$, he relied heavily on the previous work by Niebur [17] on automorphic integrals. Later, in 2000, a partial result on the missing cases in Knopp's conjecture was obtained by Wang [23] and it was resolved in 2010 by Knopp and Mawi [14], using Petersson's principal part theorem and generalised Poincaré series.

**Theorem 2** (Knopp–Mawi) *For all $r \in \mathbb{R}$, we have*

$$S_{2-r}(\Gamma, \overline{v}) \xrightarrow{\cong} H^1_{r,v}(\Gamma, \mathcal{P}).$$

A recent preprint [2] by Bruggeman and coworkers gives a similar isomorphism for a much wider class of automorphic forms. They also provide several motivations to study cocycles of real weight. One of them is a formula of Goldfeld [8] that suggests a connection between special values of derivatives of $L$-functions and cocycles. To be precise, let $f = \sum_{n \geq 1} a_n q^n$ be a Hecke cusp form of weight 2 for the group $\Gamma_0(N)$, and assume that $f$ is invariant under the Fricke involution $W_N = \begin{pmatrix} 0 & -1 \\ N & 0 \end{pmatrix}$. The $L$-function of $f$, $L_f(s)$, is defined as the analytic continuation to $\mathbb{C}$ of the Dirichlet series $\sum a_n n^{-s}$. In [2, §9.4], it is shown that Goldfeld's formula leads to the following expression:

$$-\pi i r L'_f(1) + \mathcal{O}_{r \to 0}(r^2) = \phi_{f_r}(S)(0),$$

where $f_r(z) = f(z)(\eta(z)\eta(Nz))^r$ is a cusp form of weight $2 + r$ and $S = \begin{pmatrix} 0 & 1 \\ -1 & 0 \end{pmatrix}$.

Another reason to study cocycles in the case of half-integral weight $k$ is given in [1]. Bringmann and Rolen use the non-holomorphic Eichler integral $g^*$ (this is essentially the auxiliary integral $G$ in Sect. 2 and closely connected to the cocycle $\phi_g$) to show that the function

$$Q_g\left(\frac{d}{c}\right) = L_g(\zeta_c^d; k - 1) = \sum_{n=1}^{\infty} \frac{a_g(n) \zeta_c^{nd}}{n^{k-1}},$$

defined on $\mathbb{Q}$ is a quantum modular form.

In this article, we present a new proof of Theorem 2 for positive weights $2 - r \neq 1$ that views the isomorphism in Knopp and Mawi's theorem as a duality. The key construction is a pairing between $S_{2-r}(\Gamma, \bar{v})$ and $H^1_{r,v}(\Gamma, \mathcal{P})$ which we introduce in Sect. 2. In Sect. 3 we show that this pairing is perfect, which implies Theorem 2 for the weights we consider. The proof also implies Theorem 2 for the weights $2 - r \leq 0$ and hence for all real weights, except $2 - r = 1$.

The proof proceeds as follows: Theorem 5 and Corollary 3 show that every cocycle $\phi$ in $Z^1_{r,v}(\Gamma, \mathcal{P})$ is a coboundary in $Z^1_{r,v}(\Gamma, \mathcal{Q})$, where $\mathcal{Q}$ is a larger space of functions than $\mathcal{P}$. This means that there exists $g \in \mathcal{Q}$ such that $\phi(\gamma) = g|_{r,v}\gamma - g$ for all $\gamma \in \Gamma$. If $2 - r > 0$, we assume in the next step that $\phi$ is orthogonal to all cusp forms with respect to the pairing we construct in Definition 4. Using the description of $\phi$ as a coboundary in $Z^1_{r,v}(\Gamma, \mathcal{Q})$, we use classic results from the spectral theory of automorphic forms to show that $y^{\frac{r+2}{2}} \overline{\frac{\partial g}{\partial \bar{z}}}(z)$ is in the image of the Maass weight-raising operator $K_{-r}$ (see Proposition 6). This then implies that $\phi$ is a coboundary in $Z^1_{r,v}(\Gamma, \mathcal{P})$. In the case $2 - r = 1$ only the last step of the proof fails, since some technical complications arise in the proof of Proposition 6.

One of the advantages of the new proof is that once all the constructions are in place the problem can be solved with standard techniques from the spectral theory of automorphic forms. The main references we use for spectral theory are the excellent articles [21] by W. Roelcke. Another advantage is that the proof can easily be generalised to the case of vector-valued cusp forms. We sketch this generalisation in the last section of this article.

### 1.1 Preliminaries

The group $SL_2(\mathbb{R})$ acts on $\mathbb{C} \cup \{\infty\}$ by

$$\gamma z = \frac{az+b}{cz+d}, \quad \text{for } z \in \mathbb{C},$$

and $\gamma \infty = \frac{a}{c}$, where $\gamma = \begin{pmatrix} a & b \\ c & d \end{pmatrix} \in SL_2(\mathbb{R})$. We also define the function $j(\gamma, z) = cz + d$. For any other element $\delta$ of $SL_2(\mathbb{R})$, one has the relation

$$j(\gamma\delta, z) = j(\gamma, \delta z) j(\delta, z),$$

for all $z$ in the upper half plane $\mathcal{H} = \{z = x + iy \in \mathbb{C} | y > 0\}$. Let $r \in \mathbb{R}$. Two useful functions when dealing with real weights, introduced by Petersson [18], are

$$\omega(\gamma, \delta) = \frac{1}{2\pi} \left[ -\arg(j(\gamma\delta, z)) + \arg(j(\gamma, \delta z)) + \arg(j(\delta, z)) \right]$$

and

$$\sigma_r(\gamma, \delta) = e^{2\pi i r \omega(\gamma, \delta)}.$$

Here arg, the argument, is always chosen to lie in $(-\pi, \pi]$. The value of $\omega(\gamma, \delta)$ is independent of $z$ and in $\{-1, 0, 1\}$. From the definition, it follows that

$$\sigma_r(\gamma, \delta) j(\gamma\delta, z)^r = j(\gamma, \delta z)^r j(\delta, z)^r, \quad \gamma, \delta \in \Gamma. \tag{1}$$

Here $j(\gamma, z)^r = \exp(r \cdot \log(j(\gamma, z)))$ and log is the principal branch of the complex logarithm satisfying $\log(z) = \log|z| + i \arg(z)$ for all $z \neq 0$.

In this article, we will work with modular forms with respect to a Fuchsian group of the first kind. We sketch the definition of such groups here and refer the reader to [22, §1] for a more thorough introduction. Let $\Gamma$ be a discrete subgroup of $SL_2(\mathbb{R})$ or of $SL_2(\mathbb{R})/\{\pm I\}$. A *cusp* of $\Gamma$ is any element of $\mathbb{R} \cup \{\infty\}$ that is fixed by a parabolic element of $\Gamma$, i.e. an element of $\Gamma$ that has only one fixed point in $\mathbb{R} \cup \{\infty\}$. Let $\mathcal{H}^*$ be the union of $\mathcal{H}$ with the cusps of $\Gamma$. The quotient space $\Gamma \backslash \mathcal{H}^*$ can be given the structure of a Riemann surface such that the natural projection

$$\pi : \mathcal{H} \to \Gamma \backslash \mathcal{H}^*$$

is an open map. $\Gamma$ is called a *Fuchsian group of the first kind*, if $\Gamma \backslash \mathcal{H}^*$ is compact. For the rest of this article we assume that $\Gamma \subseteq SL_2(\mathbb{R})$ is a Fuchsian group of the first kind that contains a translation. This condition is not very restrictive. Any Fuchsian group of the first kind that has cusps is conjugate to a Fuchsian group of the first kind that contains translations. For convenience, we will also assume that $\Gamma$ contains $-I$.

A *multiplier system* of weight $r$ for $\Gamma$ is a function $v : \Gamma \to \mathbb{C}$ which satisfies the consistency condition

$$v(\gamma\delta) j(\gamma\delta, z)^r = v(\gamma)v(\delta) j(\gamma, \delta z)^r j(\delta, z)^r, \quad \forall \gamma, \delta \in \Gamma$$

or equivalently

$$v(\gamma\delta) = \sigma_r(\gamma, \delta) v(\gamma) v(\delta).$$

Note that $v$ is also a multiplier system of any weight $r' \in \mathbb{R}$ with $r' \equiv r \mod 2$ and $\bar{v}$ is a multiplier system of weight $-r$. A multiplier system is called *unitary* if $|v(\gamma)| = 1$ for all $\gamma \in \Gamma$. For the rest of this article, we fix a unitary multiplier system $v$ of weight $r$.

For a function $f$ on the upper half plane $\mathcal{H}$ and $\gamma \in SL_2(\mathbb{R})$, the slash operators $|_{r,v}$ and $|_r$ are defined by

$$f|_{r,v} \gamma(z) = \bar{v}(\gamma) j(\gamma, z)^{-r} f(\gamma z)$$

and

$$f|_r \gamma(z) = j(\gamma, z)^{-r} f(\gamma z).$$

The consistency condition for $v$ implies that

$$f|_{r,v} \gamma\delta(z) = (f|_{r,v} \gamma)|_{r,v} \delta(z), \quad \forall \gamma, \delta \in \Gamma.$$

Let $q_0 = \infty$ and $q_1, \ldots, q_m$ be a set of representatives of the cusps of $\Gamma$. For every cusp $q$, the stabiliser subgroup $\Gamma_q$ is generated by $-I$ and one generator $\sigma_q \in \Gamma$. For $q = \infty$ we choose $\sigma_\infty = \begin{pmatrix} 1 & \lambda \\ 0 & 1 \end{pmatrix}$, the minimal translation matrix in $\Gamma$ with $\lambda > 0$. Let $f$ be holomorphic on $\mathcal{H}$ and invariant under $|_{r,v}$. The equation $f(z+\lambda) = v(\sigma_\infty) f(z)$ implies that $f$ has a Fourier expansion at $\infty$ of the form

$$f(z) = \sum_{n=-\infty}^{\infty} a_{n,0} \exp\left(2\pi i (n + \kappa_0) z / \lambda\right), \tag{2}$$

where $\kappa_i \in [0, 1)$ is defined for any cusp by $v(\sigma_{q_i}) = e^{2\pi i \kappa_i}$. To find the expansion at the other cusps, choose $\sigma_{q_i}$ so that if

$$A_i \sigma_{q_i} A_i^{-1} = \begin{pmatrix} 1 & \lambda_i \\ 0 & 1 \end{pmatrix},$$

where $A_i = \begin{pmatrix} 0 & -1 \\ 1 & -q_i \end{pmatrix}$, we have $\lambda_i > 0$. The Fourier expansion of $f$ at $q_i$ is then given by

$$f|_r A_i^{-1}(z) = \sum_{n=-\infty}^{\infty} a_{n,i} \exp\left(2\pi i (n + \kappa_i) z / \lambda_i\right). \tag{3}$$

**Definition 1** Let $f$ be holomorphic in $\mathcal{H}$ and invariant under $|_{r,v}$. Then $f$ is called a *modular form*[1] of weight $r$ and multiplier system $v$ with respect to $\Gamma$, if in the Fourier expansions in (2) and (3) all $a_{n,i}$ with $n + \kappa_i < 0$ are zero. If in addition all $a_{n,i}$ with $n + \kappa_i = 0$ vanish, then $f$ is called a cusp form. The set of modular forms is denoted by $M_r(\Gamma, v)$, the set of cusp forms by $S_r(\Gamma, v)$.

*Remark 1* By the main theorem of [11], the only modular form of negative weight is the zero function. By [15], the only non-zero modular forms of weight 0 are constant functions.

### 1.2 Cohomology

**Definition 2** Let $\mathcal{P}$ be the space of holomorphic functions on $\mathcal{H}$ such that there exist positive constants $K$, $A$ and $B$ with

$$|f(z)| < K(|z|^A + y^{-B}), \quad \forall z = x + iy \in \mathcal{H}.$$

A *cocycle* of weight $r$ and multiplier system $v$ with values in $\mathcal{P}$ is a function $\phi : \Gamma \to \mathcal{P}$ that satisfies

$$\phi(\gamma \delta)(z) = \phi(\gamma)|_{r,v} \delta(z) + \phi(\delta)(z), \quad \forall \gamma, \delta \in \Gamma.$$

---

[1] Another common term for modular forms that is used e.g. in [12], is *entire automorphic forms*.

We denote the space of cocycles by $Z^1_{r,v}(\Gamma, \mathcal{P})$. There is a natural map $d$ from $\mathcal{P}$ to $Z^1_{r,v}(\Gamma, \mathcal{P})$ that associates to a function $g \in \mathcal{P}$ the cocycle

$$dg : \gamma \mapsto g|_{r,v}\gamma - g.$$

A cocycle of the form $dg$ for $g \in \mathcal{P}$ is called a *coboundary* and the space of coboundaries is denoted by $B^1_{r,v}(\Gamma, \mathcal{P})$. The (first) *Eichler cohomology group* $H^1_{r,v}(\Gamma, \mathcal{P})$ is the quotient space $Z^1_{r,v}(\Gamma, \mathcal{P})/B^1_{r,v}(\Gamma, \mathcal{P})$.

A cocycle $\phi$ is called *parabolic* if for all cusps $q_i$, there exists a function $g_{q_i} \in \mathcal{P}$ such that

$$\phi(\sigma_{q_i}) = g_{q_i}|_{r,v}\sigma_{q_i} - g_{q_i}.$$

We denote the space of parabolic cocycles by $\tilde{Z}^1_{r,v}(\Gamma, \mathcal{P})$. Since coboundaries are clearly parabolic, we can form the parabolic cohomology group $\tilde{H}^1_{r,v}(\Gamma, \mathcal{P}) = \tilde{Z}^1_{r,v}(\Gamma, \mathcal{P})/B^1_{r,v}(\Gamma, \mathcal{P})$. It turns out that all cocycles are parabolic. This follows from a result that Knopp attributes to B.A. Taylor [12].

**Proposition 1** *Let $\epsilon \in \mathbb{C}$ with $|\epsilon| = 1$ and $g \in \mathcal{P}$. Then there exists an $f \in \mathcal{P}$ with*

$$\bar{\epsilon} f(z+1) - f(z) = g(z), \quad \forall z \in \mathcal{H}. \tag{4}$$

*Proof* This is Proposition 9 in [12] and a full proof is given there. We will only present the main idea here. A formal solution of (4) is given by the one-sided average

$$f(z) = -\sum_{n=0}^{\infty} \bar{\epsilon}^n g(z+n).$$

However, this sum does not always converge. Knopp uses the fact that $\mathcal{P}$ is closed under integration and differentiation to replace $g$ with a function $g' = g_1 + g_2$ such that the one-sided averages $f_1(z) = -\sum_{n=0}^{\infty} \bar{\epsilon}^n g_1(z+n)$ and $f_2(z) = -\sum_{n=0}^{\infty} \bar{\epsilon}^n g_2(z+n)$ converge and are in $\mathcal{P}$. □

**Corollary 1** *Let $\epsilon \in \mathbb{C}$ with $|\epsilon| = 1$, $s \in \mathbb{R} \setminus \{0\}$ and $g \in \mathcal{P}$. Then there exists an $f \in \mathcal{P}$ with*

$$\bar{\epsilon} f(z+s) - f(z) = g(z), \quad \forall z \in \mathcal{H}. \tag{5}$$

*Proof* First assume $s > 0$ and set $\hat{g}(z) = g(sz)$. By Proposition 1, there exists $\hat{f} \in \mathcal{P}$ that satisfies

$$\bar{\epsilon} \hat{f}(z+1) - \hat{f}(z) = \hat{g}(z), \quad \forall z \in \mathcal{H}.$$

Then $f(z) = \hat{f}(z/s)$ solves (5).

Now we treat the case $s < 0$. By the first part of this proof, there exists an $\hat{f} \in \mathcal{P}$ that satisfies

$$\epsilon \hat{f}(z - s) - \hat{f}(z) = g(z), \quad \forall z \in \mathcal{H}.$$

The function $f(z) = -\epsilon \hat{f}(z - s)$ solves (5). $\square$

**Theorem 3** ([12, p. 627]) *Every cocycle in $Z^1_{r,v}(\Gamma, \mathcal{P})$ is parabolic, i.e.*

$$Z^1_{r,v}(\Gamma, \mathcal{P}) = \tilde{Z}^1_{r,v}(\Gamma, \mathcal{P}).$$

*Proof* Let $\phi \in Z^1_{r,v}(\Gamma, \mathcal{P})$. We will show that for every parabolic $\gamma \in \Gamma$ there exists $f \in \mathcal{P}$ such that

$$\phi(\gamma) = f|_{r,v}\gamma - f. \tag{6}$$

First suppose $\gamma = \begin{pmatrix} 1 & s \\ 0 & 1 \end{pmatrix}$ is a translation by $s \in \mathbb{R} \setminus \{0\}$. Then by Corollary 1, a function $f \in \mathcal{P}$ with the desired property exists.

For the general case, let $\gamma = \begin{pmatrix} a & b \\ c & d \end{pmatrix} \in \Gamma$ and fix a cusp $q$. Then there exists an $s \in \mathbb{R} \setminus \{0\}$ such that

$$A\gamma A^{-1} = \begin{pmatrix} 1 & s \\ 0 & 1 \end{pmatrix} = U, \quad \text{where } A = \begin{pmatrix} 0 & -1 \\ 1 & -q \end{pmatrix}.$$

Replacing $z$ by $A^{-1}z$ in Eq. (6) we see that it is sufficient to show the existence of $f \in \mathcal{P}$ with

$$\overline{v(\gamma)} j(A^{-1}UA, A^{-1}z)^{-r} f(\gamma A^{-1}z) - f(A^{-1}z) = \phi(\gamma)(A^{-1}z). \tag{7}$$

Setting $\hat{f}(z) = f(A^{-1}z)$, this is equivalent to

$$\overline{v(\gamma)} j(A^{-1}UA, A^{-1}z)^{-r} \hat{f}(z + s) - \hat{f}(z) = \phi(\gamma)(A^{-1}z). \tag{8}$$

Equation (1) implies the two relations

$$1 = j(AA^{-1}U, z)^{-r} = \sigma_r(A, A^{-1}U) j(A, A^{-1}Uz)^{-r} j(A^{-1}U, z)^{-r}, \tag{9}$$
$$j(A^{-1}UA, A^{-1}z)^{-r} = \sigma_r(A^{-1}U, A) j(A^{-1}U, z)^{-r} j(A, A^{-1}z)^{-r}. \tag{10}$$

After multiplying Eq. (8) by $j(A, A^{-1}z)^r$ and using the two relations (9) and (10), we get

$$\overline{\epsilon} F(z + s) - F(z) = j(A, A^{-1}z)^r \phi(\gamma)(A^{-1}z), \tag{11}$$

where we set $F(z) = j(A, A^{-1}z)^r \hat{f}(z)$ and $\epsilon = v(\gamma)\overline{\sigma_r(A^{-1}U, A)}\sigma_r(A, A^{-1}U)$. Note that $|\epsilon| = 1$ and $j(A, A^{-1}z)^r\phi(\gamma)(A^{-1}z) \in \mathcal{P}$. The existence of such an $F \in \mathcal{P}$ again follows from Corollary 1. □

## 2 Petersson inner product

In this section we define the pairing that is essential for our proof of Theorem 2. We make use of the auxiliary integral of a cusp form of positive real weight. For weights greater than 2 it was introduced in [17] and for any positive weight it first appeared in [20], where also the transformation formula (12) is mentioned. Corollary 2 can also be deduced from results in these papers and [19] but the proof presented here is new.

**Definition 3** Let $r \in \mathbb{R}$ with $2 - r > 0$ and $g$ be a cusp form for the group $\Gamma$ of weight $2 - r$ and unitary multiplier system $\overline{v}$. The *auxiliary integral* of $g$ is defined as

$$G(z) = \left[-\int_z^\infty g(\tau)(\tau - \overline{z})^{-r} d\tau\right]^-,$$

where $[\cdot]^-$ indicates complex conjugation. The path of integration is the vertical line $p(t) = z + it$ where $t$ ranges from 0 to $\infty$.

Since $g$ decays exponentially towards $\infty$ the integral converges and $G$ is a smooth function from $\mathcal{H}$ to $\mathbb{C}$. We can define a cocycle by

$$\phi_g^\infty : \gamma \mapsto \phi_{g,\gamma}^\infty(z) = G|_{r,v}\gamma(z) - G(z).$$

**Proposition 2** *The cocycle $\phi_g^\infty$ is in $Z_{r,v}^1(\Gamma, \mathcal{P})$ and*

$$\phi_{g,\gamma}^\infty(z) = \left[\int_{\gamma^{-1}\infty}^\infty g(\tau)(\tau - \overline{z})^{-r} d\tau\right]^-, \tag{12}$$

*for all $\gamma \in \Gamma$.*

*Proof* Let $\gamma \in \Gamma$. We have

$$\overline{G(\gamma z)} = \int_\infty^{\gamma z} g(\tau)(\tau - \gamma\overline{z})^{-r} d\tau$$
$$= \int_{\gamma^{-1}\infty}^z g(\gamma\tau)(\gamma\tau - \gamma\overline{z})^{-r} d(\gamma\tau)$$
$$= j(\gamma, \overline{z})^r \int_{\gamma^{-1}\infty}^z g(\gamma\tau)j(\gamma, \tau)^{-2+r}(\tau - \overline{z})^{-r} d\tau.$$

In the last equality, we used

$$(\gamma\tau - \gamma\overline{z})^{-r} = \left(\frac{\tau - \overline{z}}{j(\gamma, \tau)j(\gamma, \overline{z})}\right)^{-r} = \frac{(\tau - \overline{z})^{-r}}{j(\gamma, \tau)^{-r} j(\gamma, \overline{z})^{-r}}.$$

To prove this, let

$$\alpha = \arg(\gamma\tau - \gamma\bar{z}) \text{ and } \beta = \arg(\tau - \bar{z}) - \arg(j(\gamma,\tau)) - \arg(j(\gamma,\bar{z})).$$

We know that $\alpha \equiv \beta \mod 2\pi$ and want to show $\alpha = \beta$. Both $(\gamma\tau - \gamma\bar{z})$ and $\tau - \bar{z}$ are in $\mathcal{H}$, so their arguments are in $(0, \pi)$. Furthermore, exactly one of $j(\gamma, \tau)$ and $j(\gamma, \bar{z})$ will be in $\mathcal{H}$ and one in $\overline{\mathcal{H}}$, so $-\pi < \beta < 2\pi$ and $0 < \alpha < -\pi$. Together with $\beta \equiv \alpha \mod 2\pi$, this implies $\alpha = \beta$. Now we use the modularity of $g$ to obtain

$$G(\gamma z) = j(\gamma, z)^r v(\gamma) \left[ \int_{\gamma^{-1}\infty}^{z} g(\tau)(\tau - \bar{z})^{-r} d\tau \right]^{-} \tag{13}$$

or $G|_{r,v}\gamma(z) = \left[ \int_{\gamma^{-1}\infty}^{z} g(\tau)(\tau - \bar{z})^{-r} d\tau \right]^{-}$. An application of Cauchy's theorem now gives us

$$\phi_{g,\gamma}^{\infty}(z) = G|_{r,v}\gamma(z) - G(z)$$
$$= \left[ \left( \int_{\gamma^{-1}\infty}^{z} - \int_{\infty}^{z} \right) g(\tau)(\tau - \bar{z})^{-r} d\tau \right]^{-}$$
$$= \left[ \int_{\gamma^{-1}\infty}^{\infty} g(\tau)(\tau - \bar{z})^{-r} d\tau \right]^{-}.$$

To see that $\phi_{g,\gamma}^{\infty}$ is in $\mathcal{P}$ first note that $(\tau - \bar{z})^{-r}$ is antiholomorphic in $\mathcal{H}$ as a function of $z$ (actually even in the slit plane $\mathbb{C} \setminus \{\mathbb{R}_{\geq 0} + \bar{\tau}\}$) and the integrals in the definition of $G$ and $\phi_g^{\infty}$ converge absolutely because $g$ is a cusp form. Therefore, $\phi_{g,\gamma}^{\infty}(z)$ is holomorphic in $\mathcal{H}$. To prove that $\phi_{g,\gamma}^{\infty}$ is in $\mathcal{P}$ one can use simple bounds for $|\tau - \bar{z}|^{-r}$. We sketch the procedure for the case $-r \geq 0$ and $\text{Im}(z) > 1$. In this case

$$|\tau - \bar{z}|^{-r} \leq |\tau - \bar{z}|^{\lceil -r \rceil} \leq \sum_{j=0}^{\lceil -r \rceil} \binom{\lceil -r \rceil}{j} |\tau|^{\lceil -r \rceil - j} |z|^j.$$

One can use this to bound $\phi_{g,\gamma}^{\infty}(z)$ by a polynomial in $|z|$. The other cases are dealt with similarly. $\square$

Let $f$ be another modular form of the weight $2 - r$ and multiplier system $\bar{v}$. Then, since $f$ is holomorphic

$$\frac{\partial Gf}{\partial \bar{z}}(z) = f(z) \frac{\partial G}{\partial \bar{z}}(z) = \overline{g(z)}(\bar{z} - z)^{-r} f(z) = (-2i)^{-r} f(z) \overline{g(z)} y^{-r}.$$

This is just a scalar times the integrand occurring in the Petersson inner product of $g$ and $f$ defined as

$$(f, g) = \int_{\Gamma \backslash H} f(z) \overline{g(z)} y^{-r} dx dy.$$

Choose a fundamental domain of $\Gamma$, $\mathcal{F}$. Then by Stokes' theorem, we have

$$(f, g) = -\frac{i}{2}\int_{\mathcal{F}} f(z)\overline{g(z)}y^{-r}\mathrm{d}\bar{z} \wedge \mathrm{d}z = C_{2-r}\int_{\partial\mathcal{F}} f(z)G(z)\mathrm{d}z,$$

for $C_{2-r} = -\frac{i}{2}(-2i)^r$. Now we choose a fundamental domain according to the following Proposition 4.2 in [3].

**Proposition 3** *The fundamental domain $\mathcal{F}$ can be chosen such that $\partial\mathcal{F}$ consists of an even number of geodesic segments $[A_i, A_{i+1}[^2$ for $i = 1, \ldots, 2n$ (the indices are taken modulo $2n$) and $\alpha_i \in \Gamma$ for $i = 1, \ldots, 2n$ such that there exists an involution of $\{1, \ldots, 2n\}$, denoted by $\tau$, such that*

1. *$\tau$ does not have any fixed points,*
2. *$\alpha_i A_i = A_{\tau(i)+1}$, $\alpha_i A_{i+1} = A_{\tau(i)}$,*
3. *$\alpha_{\tau(i)} = \alpha_i^{-1}$ and*
4. *$\alpha_i$ maps $[A_i, A_{i+1}[$ to $[A_{\tau(i)+1}, A_{\tau(i)}[$.*

*Example 1* For $\Gamma = \mathrm{SL}_2(\mathbb{Z})$, we choose the classic fundamental domain with $A_1 = \infty$, $A_2 = e^{2\pi i/3}$, $A_3 = i$, $A_4 = A_2 + 1$. Then $\alpha_1 = T = \begin{pmatrix} 1 & 1 \\ 0 & 1 \end{pmatrix}$ maps $[A_1, A_2[$ to $[A_1, A_4[$ and $\alpha_2 = S = \begin{pmatrix} 0 & 1 \\ -1 & 0 \end{pmatrix}$ maps $[A_2, A_3[$ to $[A_4, A_3[$. So $\tau$ is the permutation that swaps 1 with 4 and 2 with 3.

*Remark 2* For general Fuchsian groups $\Gamma$ of the first kind an example of such a fundamental domain is the Ford fundamental domain (see [6])

$$\mathcal{F} = \{z \in \mathcal{H} | \, |z| \leq \lambda/2 \text{ and } |j(\gamma, z)| > 1 \, \forall \gamma \in \Gamma\backslash\Gamma_\infty\}, \qquad (14)$$

where $\lambda$, the width of the cusp $\infty$, was defined in the last section. For the rest of this article, we will fix this fundamental domain for $\Gamma$.

We can restate Proposition 3 as

$$\partial\mathcal{F} = \bigsqcup_{m=1}^{n}\left([A_{i_m}, A_{i_m+1}\,[\sqcup\alpha_{i_m}]\,A_{i_m}, A_{i_m+1}]\right).$$

Thus, the Petersson inner product of $f$ and $g$ becomes

$$C_{2-r}\sum_{m=1}^{n}\left(\int_{A_{i_m}}^{A_{i_m+1}} - \int_{\alpha_{i_m} A_{i_m}}^{\alpha_{i_m} A_{i_m+1}}\right)f(z)G(z)\mathrm{d}z.$$

---

² $[A_i, A_{i+1}[$ denotes the geodesic in $\mathcal{H}$ that connects $A_i$ and $A_{i+1}$ and includes $A_i$ but not $A_{i+1}$.

Using the modularity of $f$, the second integral in the sum becomes

$$\int_{\alpha_{i_m} A_{i_m}}^{\alpha_{i_m} A_{i_m+1}} f(z)G(z)\mathrm{d}z = \int_{A_{i_m}}^{A_{i_m+1}} f(\alpha_{i_m} z)G(\alpha_{i_m} z)d(\alpha_{i_m} z)$$

$$= \int_{A_{i_m}}^{A_{i_m+1}} f(z)G|_{r,v}\alpha_{i_m}(z)\mathrm{d}z.$$

Finally, we arrive at

$$(f, g) = C_{2-r} \sum_{m=1}^{n} \int_{A_{i_m}}^{A_{i_m+1}} f(z)\left(G(z) - G|_{r,v}\alpha_{i_m}(z)\right) \mathrm{d}z$$

$$= -C_{2-r} \sum_{m=1}^{n} \int_{A_{i_m}}^{A_{i_m+1}} f(z)\phi_{g,\alpha_{i_m}}^{\infty}(z)\mathrm{d}z.$$

Motivated by the previous calculations, we define a pairing between cusp forms and cocycles:

**Definition 4** Let $2 - r > 0$, $f \in S_{2-r}(\Gamma, \bar{v})$ and $\phi \in Z_{r,v}^1(\Gamma, \mathcal{P})$. Define the pairing

$$(f, \phi) = -C_{2-r} \sum_{m=1}^{n} \int_{A_{i_m}}^{A_{i_m+1}} f(z)\phi(\alpha_{i_m})(z)\mathrm{d}z.$$

The integrals in the sum converge because $\phi(\alpha_{i_m})$ is in $\mathcal{P}$ and therefore can increase only polynomially towards the cusps, while $f$ decreases exponentially.

**Lemma 1** Let $f \in S_{2-r}(\Gamma, \bar{v})$ and $[\phi] \in H_{r,v}^1(\Gamma, \mathcal{P})$ be represented by $\phi \in Z_{r,v}^1(\Gamma, \mathcal{P})$. The value $(f, \phi)$ does not depend on a choice of representative of $[\phi]$, i.e. the pairing

$$(f, [\phi]) = (f, \phi),$$

between $S_{2-r}(\Gamma, \bar{v})$ and $H_{r,v}^1(\Gamma, \mathcal{P})$, is well defined.

*Proof* It suffices to show that if $\phi$ is a coboundary, then $(f, \phi) = 0$. If $\phi$ is a coboundary there exists a function $h \in \mathcal{P}$ with $\phi(\gamma) = h|_{r,v}\gamma - h$. We have

$$\int_{A_{i_m}}^{A_{i_m+1}} f(z)h|_{r,v}\alpha_{i_m}(z)\mathrm{d}z = \int_{A_{i_m}}^{A_{i_m+1}} f(z)j(\alpha_{i_m}, z)^{2-r}\overline{v(\alpha_{i_m})}h(\alpha_{i_m} z)d(\alpha_{i_m} z)$$

$$= \int_{A_{i_m}}^{A_{i_m+1}} f(\alpha_{i_m} z)h(\alpha_{i_m} z)d(\alpha_{i_m} z)$$

$$= \int_{\alpha_{i_m} A_i}^{\alpha_{i_m} A_{i_m+1}} f(z)h(z)\mathrm{d}z. \tag{15}$$

So

$$(f, \phi) = -C_{2-r} \sum_{m=1}^{n} \int_{A_{i_m}}^{A_{i_m+1}} f(z)\phi(\alpha_{i_m})(z)\mathrm{d}z = C_{2-r} \int_{\partial \mathcal{F}} f(z)h(z)\mathrm{d}z. \quad (16)$$

Since $(f(z)h(z))$ decays exponentially at the cusps, we can approach $(\int_{\partial \mathcal{F}} f(z)h(z)\mathrm{d}z)$ by integrals over closed paths contained in $(\mathcal{H})$. These integrals all vanish, because $(f(z)h(z))$ is holomorphic, and so $((f, \phi) = 0)$. □

**Corollary 2** *The map $f \mapsto [\phi_f^\infty]$ from $S_{2-r}(\Gamma, \overline{v})$ to $H_{r,v}^1(\Gamma, \mathcal{P})$ is injective.*

*Proof* If $[\phi_f^\infty]$ is a coboundary in $Z_{r,v}^1(\Gamma, \mathcal{P})$, then by the above calculations $0 = (f, \phi_f^\infty) = (f, f)$ and hence $f = 0$. □

## 3 The Duality theorem

In this section we prove that the pairing we defined in Lemma 1, between $S_{2-r}(\Gamma, \overline{v})$ and $H_{r,v}^1(\Gamma, \mathcal{P})$, is perfect for $0 < 2-r \neq 1$. For such weights $r$ this implies Theorem 2.

We already know that for every non-zero $f$ in $S_{2-r}(\Gamma, \overline{v})$ there exists a cocycle $\phi$ such that $(f, [\phi]) \neq 0$, since $(f, [\phi_f^\infty]) = (f, f) \neq 0$. To show that the pairing is perfect, we therefore need to prove the following theorem.

**Theorem 4** *Let $1 \neq r < 2$ and $[\phi] \in H_{r,v}^1(\Gamma, \mathcal{P})$. If $(f, [\phi]) = 0$ for all $f \in S_{2-r}(\Gamma, \overline{v})$, then $[\phi] = 0$. Together with Corollary 2 this implies that $S_{2-r}(\Gamma, \overline{v})$ and $H_{r,v}^1(\Gamma, \mathcal{P})$ are dual to each other.*

The proof of Theorem 4 will be given at the end of this section. Most constructions that follow will be valid for any real $r$ and so, if not explicitly stated otherwise, we work in this generality. In particular, we will also show Theorem 2 for $r \geq 2$.
Let $\overline{\mathcal{H}} = \mathcal{H} \cup \mathbb{R} \cup \{\infty\}$ be the closure of $\mathcal{H}$ in $\mathbb{P}^1(\mathbb{C})$. A basis of neighbourhoods of $\infty$ in $\overline{\mathcal{H}}$ is given by the sets

$$H_Y(\infty) = \{z \in \mathcal{H} | \operatorname{Im}(z) > Y\} \cup \{\infty\}.$$

Let $q$ be a cusp with $\tau_q \infty = q$ for $\tau_q \in \operatorname{SL}_2(\mathbb{R})$ such that $\tau_q^{-1} \Gamma_q \tau_q$ is generated by $T = \begin{pmatrix} 1 & 1 \\ 0 & 1 \end{pmatrix}$. Then the open sets $H_Y(q) = \tau_q H_Y(\infty)$ for $Y > 0$ form a basis of neighbourhoods of $q$.

We define a variation of the space $\mathcal{P}$ that will be useful in our proof. Let $\tilde{\mathcal{Q}}$ be the space of $C^\infty$-functions $f$ on $\mathcal{H}$ such that, for every cusp $q$ of $\Gamma$, there exists a neighbourhood $U_q \subseteq \mathcal{H}$ and $K_q, A_q, B_q > 0$ such that $f$ is holomorphic in $U_q$ and

$$|f(z)| < K_q(|z|^{A_q} + y^{-B_q}), \quad z \in U_q.$$

For the purpose of proving Theorem 4 we will actually be interested in a subspace $\mathcal{Q} \subseteq \tilde{\mathcal{Q}}$, that we introduce in Definition 5.

**Theorem 5** *Every element of $Z^1_{r,v}(\Gamma, \mathcal{P})$ is a coboundary in $Z^1_{r,v}(\Gamma, \tilde{\mathcal{Q}})$.*

*Proof* Let $\phi \in Z^1_{r,v}(\Gamma, \mathcal{P})$. We need to show that there exists a function $G \in \tilde{\mathcal{Q}}$ with $\phi(\gamma) = G|_{r,v}\gamma - G$ for all $\gamma$ in $\Gamma$. Choose $Y$ large enough, so that all the $H_Y(q)$ are disjoint and contain no elliptic fixed points. Define $U = \bigcup_{q \text{ cusp of } \Gamma} H_Y(q)$ and $V = \bigcup_{q \text{ cusp of } \Gamma} H_{2Y}(q)$. Then $U$ and $V$ are $\Gamma$-invariant. Recall that the projections $\pi(U)$ and $\pi(V)$ are open in $\Gamma \backslash \mathcal{H}^*$. By the smooth Urysohn lemma (see for example [4, Corollary 3.5.5]), there exists a smooth function $\hat{\eta}$ on $\Gamma \backslash \mathcal{H}^*$ such that $\hat{\eta}(\pi(z)) = 1$ for all $\pi(z) \in \pi(V)$ and $\hat{\eta}(\pi(z)) = 0$ for all $\pi(z)$ outside $\pi(U)$. Define $\eta(z) = \hat{\eta}(\pi(z))$ to be the pullback of $\hat{\eta}$. It is a $\Gamma$-invariant $C^\infty$-function on $\mathcal{H}$ that satisfies $\eta(z) = 1$ on $V$ and $\eta(z) = 0$ outside $U$.

We will first construct a function that has $\eta\phi$ as a coboundary. By Theorem 3, $\phi$ is a parabolic cocycle, so for every cusp $q$ there exists a function $g_q \in \mathcal{P}$ such that $\phi(\sigma_q) = g_q|_{r,v}\sigma_q - g_q$, where $\sigma_q$ is the generator of $\Gamma_q/\{\pm I\}$. We define $G$ on $U$ as follows: if $z \in H_Y(q_i)$ for some $i$ then $G(z) = g_{q_i}(z)$. If $z = \delta w$ for $\delta \in \Gamma$ and $w \in H_Y(q_i)$ we define

$$G(z) = v(\delta) j(\delta, w)^r (\phi(\delta)(w) + g_{q_i}(w)).$$

Note that this is equivalent to defining $G|_{r,v}\delta(w) = \phi(\delta)(w) + G(w)$, so once we show that the definition of $G(z)$ does not depend on the choice of $\delta$, the coboundary of $\eta G$ will be $\eta\phi$. Suppose $z = \delta w = \delta' w'$, for $\delta, \delta' \in \Gamma$ and $w, w' \in H_Y(q_i)$. We need to check that

$$v(\delta) j(\delta, w)^r (\phi(\delta)(w) + g_{q_i}(w)) = v(\delta') j(\delta', w')^r (\phi(\delta')(w') + g_{q_i}(w')).$$

Multiplying both sides by $v(\delta)^{-1} j(\delta, w)^{-r}$ and using the consistency condition of the multiplier system $v$, we see that this is equivalent to

$$\phi(\delta)(w) + g_{q_i}(w) = \left[\phi(\delta') + g_{q_i}\right]|_{r,v}(\delta'^{-1}\delta)(w).$$

This follows from the cocycle condition on $\phi$ and the choice of $g_{q_i}$. Indeed, since $w' \in \delta'^{-1}\delta H_Y(q_i) \cap H_Y(q_i) \neq \emptyset$ and since we assumed that all the $H_Y(q)$ are disjoint, $\delta'^{-1}\delta$ must fix $q_i$. Hence $\delta'^{-1}\delta = \pm\sigma_{q_i}^n$ for some $n \geq 1$. This implies

$$g_{q_i}|_{r,v}(\delta'^{-1}\delta)(w) = \phi(\delta'^{-1}\delta)(w) + g_{q_i}(w),$$

and so

$$\left[\phi(\delta') + g_{q_i}\right]|_{r,v}(\delta'^{-1}\delta)(w) = \phi(\delta)(w) - \phi(\delta'^{-1}\delta)(w) + g_{q_i}|_{r,v}(\delta'^{-1}\delta)(w)$$
$$= \phi(\delta)(w) + g_{q_i}(w).$$

So $\eta G$ is a well-defined function in $\tilde{\mathcal{Q}}$. We have thus shown that $\eta\phi$ is a coboundary in $Z^1_{r,v}(\Gamma, \tilde{\mathcal{Q}})$.

It remains to show that $(1 - \eta)\phi$ is a coboundary. We first construct a partition of unity on $\mathcal{H}$ that is $\Gamma$-invariant. The construction we describe here is due to Gunning [9]. Since $\Gamma$ acts discontinuously on $\mathcal{H}$, every $z \in \mathcal{H}$ has a neighbourhood $O_z$ such that $\gamma O_z = O_z$ if $\gamma \in \Gamma_z$ (the stabiliser of $z$) and $\gamma O_z \cap O_z = \emptyset$ if $\gamma \in \Gamma \setminus \Gamma_z$. Let $V$ be as in the construction of $\eta$, a $\Gamma$-invariant open set that contains all cusps of $\Gamma$ with $\eta|_V = 1$. Since $\Gamma \backslash \mathcal{H}^*$ is compact, there exist $z_1, \ldots, z_n \in \mathcal{H}$ such that the sets $\pi(O_{z_i})$ together with $\pi(V)$ cover $\Gamma \backslash \mathcal{H}^*$. Let $\hat{\epsilon}_1, \ldots, \hat{\epsilon}_n, \hat{\epsilon}_V$ be a partition of unity corresponding to this cover, i.e. smooth functions supported in $\pi(O_{z_1}), \ldots, \pi(O_{z_n})$ and $\pi(V)$, respectively, satisfying

$$\sum_{i=1}^{n} \hat{\epsilon}_i(\pi(z)) + \hat{\epsilon}_V(\pi(z)) = 1, \quad \forall z \in \mathcal{H}.$$

We define functions $H_1, \ldots, H_n$ on $\mathcal{H}$ as follows. If there exists $g_i(z) \in \Gamma$ such that $g_i(z)z \in O_{z_i}$ we set

$$H_i(z) = -(1 - \eta(z)) \frac{\epsilon_i(z)}{|\Gamma_i|} \sum_{g \in \Gamma_i} \phi(g \cdot g_i(z))(z),$$

where $(\Gamma_i)$ is the stabiliser of $(z_i)$, $(|\Gamma_i|)$ is its order and $(\epsilon_i(z) = \hat{\epsilon}_i(\pi(z)))$. This does not depend on the choice of $g_i(z)$: if $\gamma z \in O_{z_i}$ with $\gamma \in \Gamma$, then we must have $\gamma^{-1} g_i(z) \in \Gamma_i$. Thus, the set $\Gamma_i g_i(z)$ is equal to $\Gamma_i \gamma$ and we see that a different choice of $g_i(z)$ just permutes the summands in the definition of $H_i(z)$. If no such $g_i(z) \in \Gamma$ exists, we set $H_i(z) = 0$.

Clearly, $H_i$ is a function in $\tilde{\mathcal{Q}}$ and defining $H = \sum_{i=1}^{n} H_i$, we will see that $H|_{r,v} \gamma(z) - H(z) = (1 - \eta(z))\phi(\gamma)(z)$ for all $\gamma \in \Gamma$ and $z \in \mathcal{H}$. First note that if $z$ is in $V$, then $H(z)$ and $H(\gamma z)$ vanish and so does $(1 - \eta(z))\phi(\gamma)(z)$. If $z$ is not in $V$, we have

$$H|_{r,v} \gamma(z) = -(1 - \eta(z)) \sum_i \frac{\epsilon_i(\gamma z)}{|\Gamma_i|} \sum_{g \in \Gamma_i} \phi(g \cdot g_i(\gamma z))|_{r,v}\gamma(z),$$

where the first sum is over all $i$ such that there exists a $g_i(\gamma z) \in \Gamma$ with $g_i(\gamma z)\gamma z \in O_{z_i}$. Now we choose $g_i(\gamma z) = g_i(z)\gamma^{-1}$, to get that $H|_{r,v}\gamma(z)$ equals

$$= -(1 - \eta(z)) \sum_i \frac{\epsilon_i(z)}{|\Gamma_i|} \sum_{g \in \Gamma_i} [\phi(g \cdot g_i(z))(z) - \phi(\gamma)(z)]$$
$$= (1 - \eta(z))(\phi(\gamma)(z) + H(z)).$$

□

In the definition of $\tilde{\mathcal{Q}}$, the constants $K_q, A_q, B_q$ may vary from cusp to cusp; in the following definition, we impose stricter growth conditions, requiring the constants to be fixed.

**Definition 5** Let $\mathcal{Q}$ be the space of functions $F$ in $\tilde{\mathcal{Q}}$ such that there exist positive constants $K$, $A$, $B$ with

$$|F(z)| < K(|z|^A + y^{-B}), \quad \forall z \in \mathcal{H}.$$

Note that the functions in $\mathcal{P}$ are the holomorphic functions in $\mathcal{Q}$.

**Proposition 4** *Let $F$ be in $\tilde{\mathcal{Q}}$. If $\gamma \mapsto F|_{r,v}\gamma - F = \psi(\gamma)$ is in $Z^1_{r,v}(\Gamma, \mathcal{P})$ then $F$ is in $\mathcal{Q}$.*

*Proof* This proof is similar to the proof of the main theorem of [13]. Let $M$ be the set of matrices $\gamma$ in $\Gamma$ with $\lambda/2 \leq \text{Re}(\gamma i) < \lambda/2$. $M$ is a complete set of coset representatives of $\Gamma_\infty \setminus \Gamma$. We need a technical lemma from [12]:

**Lemma 2** (Lemma 8 in [12]) *There exist positive constants $K_1$, $A_1$, $B_1$ such that for all $\tau \in \overline{\mathcal{F}} \cap \mathcal{H}$ and all $\gamma \in M$*

$$|\psi(\gamma)(\tau)| < K_1(Im(\gamma\tau)^{A_1} + Im(\gamma\tau)^{-B_1}).$$

Since only finitely many cusps are in $\overline{\mathcal{F}}$ and since the real part of $z \in \mathcal{F}$ is bounded, we can also find positive $K_2$, $A_2$, $B_2$ with

$$|F(\tau)| < K_2(Im(\tau)^{A_2} + Im(\tau)^{-B_2}), \quad \forall \tau \in \overline{\mathcal{F}} \cap \mathcal{H}. \tag{17}$$

As in the proof of Theorem 5, we use the fact that $\psi$ is parabolic and hence there exists a function $g_\infty \in \mathcal{P}$ such that $\psi(\sigma_\infty) = g_\infty|_{r,v}\sigma_\infty - g_\infty$. The equation $F|_{r,v}\sigma_\infty - F = \psi(\sigma_\infty)$ implies

$$(F - g_\infty)|_{r,v}\sigma_\infty - (F - g_\infty) = 0.$$

$F$ is in $\mathcal{Q}$ if and only if $F - g_\infty$ is in $\mathcal{P}$, so we can assume without loss of generality that $F(z+\lambda) = v(\sigma_\infty)F(z)$. Let $z \in \mathcal{H}$. There exists $\tau \in \overline{\mathcal{F}}$ and $\gamma \in \Gamma$ such that $z = \gamma\tau$. Since $M$ is a complete set of representatives of $\Gamma_\infty \setminus \Gamma$, there is an integer $m$ and $\delta \in M$ such that $z = \sigma_\infty^m \delta\tau$. If $\delta = I$ then we can deduce

$$|F(z)| < K_2(Im(\tau)^{A_2} + Im(z)^{-B_2}),$$

from Eq. (17) and the fact that $|F|$ is $\Gamma_\infty$-invariant. Suppose $\delta = \begin{pmatrix} a & b \\ c & d \end{pmatrix}$ is not the identity. Then $c \neq 0$, because the only member of $M$ that fixes $\infty$ is $I$. We have

$$|F(z)| = |F(\sigma_\infty^m \delta\tau)| = |F(\delta\tau)| \tag{18}$$

$$\leq |j(\delta,\tau)|^r (|F(\tau)| + |\psi(\delta)(\tau)|) \tag{19}$$

$$< |j(\delta,\tau)|^r [K_2(Im(\tau)^{A_2} + Im(\tau)^{-B_2}) + K_1(Im(\delta\tau)^{A_1} + Im(\delta\tau)^{-B_1})]. \tag{20}$$

By our choice of fundamental domain we have $|j(\delta, \tau)| \geq 1$, since $\delta \notin \Gamma_\infty$. So $y = \text{Im}(z) = \frac{\text{Im}(\tau)}{|j(\delta,\tau)|^2} \leq \text{Im}(\tau)$. On the other hand, using $\tau = \delta^{-1}\sigma_\infty^{-m}z$, we have $\text{Im}(\tau) = \frac{y}{|j(\delta^{-1}\sigma_\infty^{-m},z)|^2}$ and

$$\left|j(\delta^{-1}\sigma_\infty^{-m}, z)\right|^2 = |-cz + cm\lambda + a|^2 = c^2 y^2 + (cm\lambda + a - cx)^2 \geq cy^2 > c_0 y^2,$$

where $c_0 > 0$ depends only on $\Gamma$. Such a $c_0$ exists because $\Gamma$ is discrete. Therefore, $y \leq \text{Im}(\tau) < c_0^{-1} y^{-1}$, $\text{Im}(\tau)^{A_2} < c_0^{-A_2} y^{-A_2}$ and $\text{Im}(\tau)^{-B_2} \leq y^{B_2}$. Also $|j(\delta, \tau)|^r = \left(\frac{y}{\text{Im}(\tau)}\right)^{-r/2}$ is either $\leq 1$ (if $r \leq 0$) or $\leq c_0^{-r/2} y^{-r}$ (if $r \geq 0$). These inequalities inserted into (20) lead to the desired inequality of the form

$$|F(z)| < K(|z|^A + y^{-B}),$$

for positive constants $K, A, B$ and all $z \in \mathcal{H}$. □

**Corollary 3** *Every cocycle in $Z^1_{r,v}(\Gamma, \mathcal{P})$ is a coboundary in $Z^1_{r,v}(\Gamma, \mathcal{Q})$.*

Let $\phi \in Z^1_{r,v}(\Gamma, \mathcal{P})$. By Corollary 3 there exists a function $g \in \mathcal{Q}$ such that $g|_{r,v}\gamma - g = \phi(\gamma)$ for all $\gamma \in \Gamma$. By the same calculation as in Eq. (16), for any $f \in S_{2-r}(\Gamma, \overline{v})$, we have

$$(f, \phi) = -C_{2-r} \sum_{m=1}^{n} \int_{A_{i_m}}^{A_{i_m+1}} f(z)(g|_{r,v}\alpha_{i_m}(z) - g(z))dz$$

$$= C_{2-r} \int_{\partial \mathcal{F}} f(z)g(z)dz$$

$$= C_{2-r} \int_{\mathcal{F}} \frac{\partial g}{\partial \bar{z}} d\bar{z} \wedge f(z)dz.$$

Here we note again that the integrals above exist because $g$ can only increase polynomially towards the cusps of $\Gamma$, while $f$ decreases exponentially.

### 3.1 Spectral theory of automorphic forms

To carry out the proof of Theorem 4, we will apply spectral theory. We only give a very brief introduction here; for more details and proofs, see the exposition [21] by Roelcke. In these articles Roelcke uses a variation of the slash operator which we denote by $|_{r,v}^R$

$$f|_{r,v}^R \gamma(z) = \left(\frac{j(\gamma, \bar{z})}{j(\gamma, z)}\right)^{r/2} \overline{v}(\gamma) f(\gamma z).$$

The connection to our slash operator is given by the following lemma:

**Lemma 3** *Let $f : \mathcal{H} \to \mathbb{C}$, $F(z) = y^{\frac{r}{2}} f(z)$ and $\gamma \in \Gamma$. Then*

$$y^{\frac{r}{2}} \left( f|_{r,v} \gamma(z) \right) = F|^R_{r,v} \gamma(z).$$

*So a function $f$ is invariant under $|_{r,v}$ if and only if $F(z) = y^{\frac{r}{2}} f(z)$ is invariant under $|^R_{r,v}$.*

**Definition 6** Let $H_{r,v} = H_r(\Gamma, v)$ be the Hilbert space of functions $f$ that are invariant under $|^R_{r,v}$ and have finite norm with respect to the scalar product

$$(f_1, f_2)^R = \int_{\mathcal{F}} f_1(z) \overline{f_2(z)} \frac{dx\,dy}{y^2}.$$

The weight $r$ hyperbolic Laplacian and the Maass weight-raising and weight-lowering operators are defined as

$$\Delta_r = -(z - \bar{z})^2 \frac{\partial^2}{\partial z \partial \bar{z}} - \frac{r}{2}(z - \bar{z}) \left( \frac{\partial}{\partial z} + \frac{\partial}{\partial \bar{z}} \right),$$

$$K_r = (z - \bar{z}) \frac{\partial}{\partial z} + \frac{r}{2},$$

$$\Lambda_r = (z - \bar{z}) \frac{\partial}{\partial \bar{z}} + \frac{r}{2}.$$

Before we sum up the main properties of these operators in Proposition 5, we recall some definitions from operator theory.

**Definition 7** Let $H$ and $H'$ be Hilbert spaces and let $T$ be a linear operator from a subspace $D$ of $H$ to $H'$. $T$ is called *closed* if, for every sequence $x_n$ in $D$ that converges to $x \in H$ such that $T x_n$ converges to $y \in H'$, we have $x \in D$ and $Tx = y$.

**Definition 8** If $D$ is dense in $H$ then for any operator $T$ from $D$ to $H$, we can define its *adjoint* $T^*$ on the domain

$$\{y \in H : x \mapsto \langle Tx, y \rangle \text{ is continuous on } D\}.$$

Any $y$ in this set defines a linear functional on $D$ by $\phi_y : x \mapsto \langle Tx, y \rangle$. This functional can be extended to $H$ and by the Riesz representation theorem there exists $z \in H$ such that $\phi_y(x) = \langle x, z \rangle$ for all $x$ in $H$. We define $T^* y = z$.

An operator is called *self-adjoint* if it is equal to its adjoint. An operator is called *essentially self-adjoint* if $T \subseteq T^* = (T^*)^*$, where $T \subseteq T^*$ means that $T^*$ extends $T$.

Let

$$\mathcal{D}^2_r = \{f \in H_{r,v} | f \text{ twice differentiable and } -\Delta_r f \in H_{r,v}\}.$$

**Proposition 5** *(i) $\Delta_r : \mathcal{D}^2_r \to H_{r,v}$ is essentially self-adjoint. It has a self-adjoint extension to a dense subset of $H_{r,v}$ that we denote by $\tilde{\mathcal{D}}_r$.*

(ii) *The eigenfunctions of $\Delta_r$ are smooth (in fact they are real analytic).*
(iii) $K_r : \mathcal{D}_r^2 \to H_{r+2,v}$ *and* $\Lambda_r : \mathcal{D}_r^2 \to H_{r-2,v}$ *can be extended to closed operators defined on $\tilde{\mathcal{D}}_r$. For $f \in \tilde{\mathcal{D}}_r$ and $g \in \tilde{\mathcal{D}}_{2+r}$, we have*

$$(K_r f, g)^R = (f, \Lambda_{2+r} g)^R.$$

(iv)
$$-\Delta_r = \Lambda_{r+2} K_r - \frac{r}{2}\left(1 + \frac{r}{2}\right) = K_{r-2} \Lambda_r + \frac{r}{2}\left(1 - \frac{r}{2}\right).$$

*Proof* For proofs of the statements (i), (iii) and (iv) see [21]. (i) is Satz 3.2, (iii) follows from the discussion after the proof of Lemma 6.2 on page 332 and (iv) is Eq. (3.4) on page 305. Statement (ii) follows from the fact that $\Delta_r$ is an elliptic operator and elliptic regularity applies. For an introduction to the theory of elliptic operators, see [7]. The result needed here is Corollary 8.11 in [7]. □

**Definition 9** A *cuspidal Maass wave form* in $H_{r,v}$ with eigenvalue $\lambda$ is an eigenfunction of $-\Delta_r$ with eigenvalue $\lambda$ that decays exponentially at the cusps of $\Gamma$.

*Remark 3* By [21, Satz 5.2], all eigenfunctions in $H_{r,v}$ of $-\Delta_r$ of eigenvalue $\frac{r}{2}(1-\frac{r}{2})$ are of the form $y^{\frac{r}{2}} f$, where $f$ is a modular form in $M_r(\Gamma, v)$ that has finite Petersson norm, i.e. $(f, f) < \infty$. If $f$ is a cusp form, then $y^{\frac{r}{2}} f$ is a cuspidal Maass wave form.

The main result in [21] is a spectral decomposition of $\Delta_r$. For this purpose we introduce the Eisenstein series. Let $q$ be a cusp of $\Gamma$, $\sigma_q$ the generator of $\Gamma_q/\{\pm I\}$ and $A_q \in \mathrm{SL}_2(\mathbb{R})$ chosen such that $q = A_q^{-1}\infty$. The cusp $q$ is called *singular for the multiplier system* $v$, if $v(\sigma_q) = 1$ and *regular for $v$* otherwise. Let $q_1, \ldots, q_{m^*}$ be a set of representatives of the cusps of $\Gamma$ that are singular for $v$. For each of these cusps, we define the Eisenstein series

$$E_{r,v}^q(z, s) = \frac{1}{2} \sum_{M \in \Gamma_q \backslash \Gamma} \sigma_r(A_q, M)^{-1} \overline{v(M)} \left(\frac{j(A_q M, \bar{z})}{j(A_q M, z)}\right)^{r/2} (\mathrm{Im}\, A_q M z)^s.$$

The definition of $E_{r,v}^q$ depends on the choice of $A_q$, but a different choice of $A_q$ will only multiply the Eisenstein series by a constant of absolute value 1. The series above converge absolutely and uniformly for $(z, s)$ in sets of the form $K \times \{s | \mathrm{Re}\, s \geq 1 + \epsilon\}$, where $K$ is a compact subset of $\mathbb{C}$ and $\epsilon > 0$. For a fixed $s$ with real part $\geq 1 + \epsilon$, one can use the absolute and uniform convergence of the series to see that $E_{r,v}^q(\cdot, s)$ is invariant under $|_{r,v}^R$ and that

$$-\Delta_r E_{r,v}^q(\cdot, s) = s(1-s) E_{r,v}^q(\cdot, s).$$

These series can be meromorphically continued and play an important role in the spectral decomposition of $\Delta_r$.

**Theorem 6** (i) *For fixed $z \in \mathcal{H}$, the Eisenstein series $E_{r,v}^q(z, \cdot)$ can be meromorphically continued to the whole complex plane.*
(ii) *If, for one fixed $z$, $E_{r,v}^q(z, \cdot)$ has a pole of order $n$ at $s_0$, then the function $f(z) := \lim_{s \to s_0} (s - s_0)^n E_{r,v}^q(z, s)$ is real analytic, invariant under $|_{r,v}^R$ and satisfies*

$$-\Delta_r f = s_0(1 - s_0) f.$$

*If $n$ is chosen so that $f(z)$ has no poles in $\mathcal{H}$, then $f$ grows at most polynomially at each cusp of $\Gamma$, i.e. if $q$ is a cusp of $\Gamma$ and $\tau_q \infty = q$ for $\tau_q \in \mathrm{SL}_2(\mathbb{R})$, then there exists $A \in \mathbb{R}$ such that $f|_r \tau_q(z) = \mathcal{O}(y^A)$ as $y \to \infty$.*
*In particular, if $E_{r,v}^q(z, s)$ is holomorphic at $s = s_0$, then*

$$-\Delta_r E_{r,v}^q(\cdot, s_0) = s_0(1 - s_0) E_{r,v}^q(\cdot, s_0).$$

*Furthermore, we have the following equalities:*

$$K_r E_{r,v}^q(\cdot, s_0) = \left(\frac{r}{2} + s_0\right) E_{r+2,v}^q(\cdot, s_0), \tag{21}$$

$$\Lambda_r E_{r,v}^q(\cdot, s_0) = \left(\frac{r}{2} - s_0\right) E_{r-2,v}^q(\cdot, s_0). \tag{22}$$

*The poles of $E_{r,v}^q(z, \cdot)$ in the half plane defined by $\mathrm{Re}\, s \geq \frac{1}{2}$ are all simple and in the interval $(\frac{1}{2}, 1]$. In particular, there are no poles on the line $\mathrm{Re}\, s = \frac{1}{2}$.*

**Theorem 7** (Spectral expansion) *Let $f \in \tilde{\mathcal{D}}_r$ and $e_n$ be a maximal orthonormal system of eigenfunctions[3] of $\Delta_r$. Then $f$ has a spectral expansion*

$$f = \sum_n (e_n, f)^R e_n + \sum_{i=1}^{m^*} \frac{1}{4\pi} \int_{-\infty}^{\infty} \left(E_{r,v}^{q_i}\left(\cdot, \frac{1}{2} + i\rho\right), f\right)^R E_{r,v}^{q_i}\left(z, \frac{1}{2} + i\rho\right) d\rho.$$

*If $f$ has compact support mod $\Gamma$, i.e. $\pi(\mathrm{supp}(f))$ is compact in $\Gamma \backslash \mathcal{H}^*$, then both parts of the spectral expansion, $\sum (e_n, f)^R e_n$ and $\sum_{i=1}^{m^*} \frac{1}{4\pi} \int_{-\infty}^{\infty} (E_{r,v}^{q_i}(\cdot, \frac{1}{2} + i\rho), f)^R E_{r,v}^{q_i}(z, \frac{1}{2} + i\rho) d\rho$, converge absolutely and uniformly on compact subsets of $\mathcal{H}$.*

Both the properties of Eisenstein series and the spectral expansion are proved in the second part of [21]. The Theorem we state is a combination of Satz 7.2 and the second part in Satz 12.3.

We turn back to the proof of Theorem 4: Let $[\phi] \in H_{r,v}^1(\Gamma, \mathcal{P})$ be represented by $\phi \in Z_{r,v}^1(\Gamma, \mathcal{P})$. By Corollary 3, there exists a function $g \in \mathcal{Q}$ such that

$$\phi(\gamma) = g|_{r,v}\gamma - g, \quad \forall \gamma \in \Gamma. \tag{23}$$

---

[3] An *orthonormal system* of eigenfunctions of an operator $T$ on a Hilbert space $H$ is a set of eigenfunctions of $T$ that are pairwise orthogonal and have norm 1.

By applying $\frac{\partial}{\partial \bar{z}}$ to (23), we see that

$$\frac{\partial g}{\partial \bar{z}}(z) = \overline{v(\gamma)} j(\gamma, z)^{-r} j(\gamma, \bar{z})^{-2} \frac{\partial g}{\partial \bar{z}}(\gamma z).$$

A short calculation shows that the function

$$G : z \mapsto y^{\frac{r+2}{2}} \overline{\frac{\partial g}{\partial \bar{z}}(z)} \tag{24}$$

is invariant under $|_{2-r,\bar{v}}^{R}$. Moreover, $G$ vanishes in a neighbourhood of every cusp since $g$ is holomorphic there, so $G$ has compact support mod $\Gamma$ and is in $H_{2-r,\bar{v}}$.

To prove Theorem 4, we have to show that if $\phi$ is orthogonal to $S_{2-r}(\Gamma, \bar{v})$, then $g \in \mathcal{Q}$ can be chosen to be holomorphic. This implies that $\phi$ is a coboundary in $Z_{r,v}^{1}(\Gamma, \mathcal{P})$.

**Lemma 4** *Let $2 - r > 0$ and $\phi$, $g$ and $G$ be as above. Then $(f, \phi) = 0$ for all $f \in S_{2-r}(\Gamma, \bar{v})$ if and only if $(\tilde{f}, G)^{R} = 0$ for all cuspidal Maass wave forms $\tilde{f}$ with eigenvalue $\frac{r}{2}(1 - \frac{r}{2})$.*

*Proof* We have the equality

$$\frac{i}{2C_{2-r}}(f, \phi) = \frac{i}{2} \int_{\mathcal{F}} \bar{\partial} g \wedge f(z) \mathrm{d}z = \int_{\mathcal{F}} y^{\frac{2-r}{2}} f(z) \overline{G(z)} \frac{\mathrm{d}x \mathrm{d}y}{y^{2}} = (y^{\frac{2-r}{2}} f, G)^{R},$$

so $(f, \phi) = 0$ for all $f \in S_{2-r}(\Gamma, \bar{v})$ if and only if $(\tilde{f}, G)^{R} = 0$ for all functions $\tilde{f}$ of the form $y^{\frac{2-r}{2}} f$, $f \in S_{2-r}(\Gamma, \bar{v})$. According to Remark 3, these functions are exactly the cuspidal Maass wave forms of eigenvalue $\frac{r}{2}(1 - \frac{r}{2})$. □

We can now use spectral theory to characterise functions which are orthogonal to cuspidal Maass wave forms of eigenvalue $\frac{r}{2}(1 - \frac{r}{2})$.

**Proposition 6** *Let $2 - r \neq 1$ and $H$ be a smooth function in $H_{2-r,\bar{v}}$ with compact support mod $\Gamma$. Then the following are equivalent:*

(i) *$(\tilde{f}, H)^{R} = 0$ for all cuspidal Maass wave forms $\tilde{f}$ with eigenvalue $\frac{r}{2}(1 - \frac{r}{2})$.*
(ii) *$H = K_{-r} F + K_{-r} E$, where $F$ is a smooth function in $H_{-r,\bar{v}}$ and $E$ is a linear combination of the functions $E_{-r,\bar{v}}^{qi}(z, \frac{r}{2})$.*
*If $2 - r > 1$ or $2 - r < 0$ this implies $E = 0$.*

*Remark 4* By [11,15] we have $S_{2-r}(\Gamma, \bar{v}) = \{0\}$, if $2 - r \leq 0$. Since, by [21, Satz 5.2], all cuspidal Maass wave forms of eigenvalue $\frac{r}{2}(1 - \frac{r}{2})$ are of the form $y^{\frac{r}{2}} f$, where $f \in S_{2-r}(\Gamma, \bar{v})$, the first condition is always satisfied in the case $2 - r \leq 0$.

*Proof* (i)⇒(ii): By [21, Satz 6.3] there is a maximal orthonormal system of eigenfunctions of $\Delta_{2-r}$ consisting of the following:

1. Images of eigenfunctions of $\Delta_{-r}$ under the Maass raising operator $K_{-r} = (z - \bar{z})\frac{\partial}{\partial z} - \frac{r}{2}$. We denote these by $K_{-r} e_{n}$. By [21, Satz 6.3] these eigenfunctions cannot have eigenvalue $\frac{r}{2}(1 - \frac{r}{2})$.

2. A (finite) orthonormal basis of the eigenfunctions of eigenvalue $\frac{r}{2}(1 - \frac{r}{2})$. By Remark 3 this set is of the form $\{y^{\frac{2-r}{2}} f_1, \ldots, y^{\frac{2-r}{2}} f_N\}$, where the $f_i$ form an orthonormal basis of the subspace of $M_{2-r}(\Gamma, \overline{v})$ of modular forms with finite Petersson norm. If $2 - r \geq 1$ this subspace is equal to $S_{2-r}(\Gamma, \overline{v})$, while for $2 - r < 1$ every modular form in $M_{2-r}(\Gamma, \overline{v})$ has finite Petersson norm.

Hence by Theorem 7, the spectral expansion of $H$ is of the form

$$H = \underbrace{\sum_n (K_{-r} e_n, H)^R K_{-r} e_n}_{=K_{-r} F_1} + \underbrace{\sum_{i=1}^N (y^{\frac{2-r}{2}} f_i, H)^R y^{\frac{2-r}{2}} f_i}_{=y^{\frac{2-r}{2}} \tilde{E}}$$

$$+ \underbrace{\sum_{i=1}^{m^*} \frac{1}{4\pi} \int_{-\infty}^{\infty} \left(E^{q_i}_{2-r,\overline{v}}\left(\cdot, \frac{1}{2} + i\rho\right), H\right)^R E^{q_i}_{2-r,\overline{v}}\left(z, \frac{1}{2} + i\rho\right) d\rho}_{=\tilde{F}_2}.$$

Here we used that $\sum_n (K_{-r} e_n, H)^R K_{-r} e_n$ converges absolutely and uniformly on compacta to swap differentiation and summation and write it as $K_{-r} F_1 = K_{-r} \left(\sum_n (K_{-r} e_n, H)^R e_n\right)$.

We now show that $\tilde{F}_2 = K_{-r} F_2$ for a smooth function $F_2 \in H_{-r, \overline{v}}$: Applying Eq. (21) twice and using Proposition 5, we see

$$\int_{-\infty}^{\infty} \left(E^{q_i}_{2-r,\overline{v}}\left(\cdot, \frac{1}{2} + i\rho\right), H\right)^R E^{q_i}_{2-r,\overline{v}}\left(z, \frac{1}{2} + i\rho\right) d\rho$$
$$= \int_{-\infty}^{\infty} \left(\frac{1-r}{2} + i\rho\right)^{-2} \underbrace{\left(K_{-r} E^{q_i}_{-r,\overline{v}}\left(\cdot, \frac{1}{2} + i\rho\right), H\right)^R}_{=(E^{q_i}_{-r,\overline{v}}, \Lambda_{2-r} H)^R} K_{-r} E^{q_i}_{-r,\overline{v}}\left(z, \frac{1}{2} + i\rho\right) d\rho.$$

If $r \neq 1$

$$F_2^i(z) = \int_{-\infty}^{\infty} \left(\frac{1-r}{2} + i\rho\right)^{-2} (E^{q_i}_{-r,\overline{v}}, \Lambda_{2-r} H)^R E^{q_i}_{-r,\overline{v}}\left(z, \frac{1}{2} + i\rho\right) d\rho \quad (25)$$

converges absolutely and uniformly on compacta. To see this note the integrand can be bounded above by

$$\left|\frac{1-r}{2}\right|^{-2} \cdot \left|(E^{q_i}_{-r,\overline{v}}, \Lambda_{2-r} H)^R E^{q_i}_{-r,\overline{v}}\left(z, \frac{1}{2} + i\rho\right)\right|,$$

and

$$\int_{-\infty}^{\infty} (E^{q_i}_{-r,\overline{v}}, \Lambda_{2-r} H)^R E^{q_i}_{-r,\overline{v}}\left(z, \frac{1}{2} + i\rho\right) d\rho,$$

converges absolutely and uniformly on compacta as it occurs in the spectral expansion of $\Lambda_{2-r}H$. So when we apply $K_{-r}$ to $F_2 = \sum_{i=1}^{m^*} \frac{1}{4\pi} F_2^i$, we can swap it with the integral and obtain

$$K_{-r} F_2 = \tilde{F}_2.$$

$F_2$ is clearly in $H_{-r,\bar{v}}$ by the bound we used for the $F_2^i$. We have thus shown that

$$H = K_{-r}F + y^{\frac{2-r}{2}} \tilde{E}, \text{ where } F = F_1 + F_2 \in H_{-r,\bar{v}}. \tag{26}$$

To see that $F$ is smooth we apply $\Lambda_{2-r}$ to (26) and obtain

$$\Lambda_{2-r} H = \Lambda_{2-r} K_{-r} F + \Lambda_{2-r}(y^{\frac{2-r}{2}} \tilde{E}) = -\Delta_{-r} F - \frac{r}{2}(1 - \frac{r}{2}) F + \Lambda_{2-r}(y^{\frac{2-r}{2}} \tilde{E}).$$

We see that $F$ is a solution of an elliptic differential equation and so, by elliptic regularity, $F$ is smooth.

It remains to show that $y^{\frac{2-r}{2}} \tilde{E}$ is in the image of $K_{-r}$. Since $H$ is orthogonal to all cuspidal Maass wave forms with eigenvalue $\frac{r}{2}(1 - \frac{r}{2})$, we see that in the expansion

$$\tilde{E} = \sum_{i=1}^{N} (y^{\frac{2-r}{2}} f_i, H)^R f_i$$

only the $f_i \in M_{2-r}(\Gamma, \bar{v})$ that are orthogonal to $S_{2-r}(\Gamma, \bar{v})$ can occur. Hence $\tilde{E}$ must be orthogonal to $S_{2-r}(\Gamma, \bar{v})$ and has finite Petersson norm. If $2 - r \geq 1$ this implies $\tilde{E} = 0$. If $2 - r < 0$, we have $M_{2-r}(\Gamma, \bar{v}) = \{0\}$ by [11], so in this case we also have $\tilde{E} = 0$. We are left with the case $0 \leq 2 - r < 1$. In this case all modular forms in $M_{2-r}(\Gamma, \bar{v})$ have finite Petersson norm, so $\tilde{E}$ can be any form in the orthogonal complement of $S_{2-r}(\Gamma, \bar{v})$. We can appeal to [21, Satz 11.2], to see that $\tilde{E}$ is a linear combination of residues of Eisenstein series at $s = \frac{r}{2}$. Therefore, there exist $a_i \in \mathbb{C}$ with

$$y^{\frac{2-r}{2}} \tilde{E}(z) = \sum_{i=1}^{m^*} a_i \operatorname{Res}_{s=\frac{r}{2}}(E_{2-r,\bar{v}}^{q_i}(z, s)).$$

Note that we can restrict the sum on the right-hand side to include only Eisenstein series that have a pole at $s = \frac{r}{2}$. On the other hand, Eisenstein series of weight $-r$ never have a pole at $s = \frac{r}{2}$ by [21, Satz 13.2], since $-r < -1$. Equation (21) now implies

$$\operatorname{Res}_{s=\frac{r}{2}}(E_{2-r,\bar{v}}^{q_i}(z, s)) = \lim_{s \to \frac{r}{2}} (s - \frac{r}{2}) E_{2-r,\bar{v}}^{q_i}(z, s) \tag{27}$$

$$= \lim_{s \to \frac{r}{2}} K_{-r} E_{-r,\bar{v}}^{q_i}(z, s) = K_{-r} E_{-r,\bar{v}}^{q_i}(z, \frac{r}{2}). \tag{28}$$

Setting $E = \sum_{i=1}^{m^*} a_i E_{-r,\overline{v}}^{q_i}(z, \frac{r}{2})$ we can confirm statement (ii).

(ii)$\Rightarrow$(i): Let $H = K_{-r}F + K_{-r}E$ as described in (ii) and let $\tilde{f}$ be a cuspidal Maass wave form with eigenvalue $\frac{r}{2}(1 - \frac{r}{2})$. From the first part of the proof we know that $K_{-r}E$ has the form $y^{\frac{2-r}{2}}\tilde{E}$, where $\tilde{E} \in M_{2-r}(\Gamma, \overline{v})$ is orthogonal to $S_{2-r}(\Gamma, \overline{v})$. This implies that $y^{\frac{2-r}{2}}\tilde{E}$ is orthogonal to $\tilde{f}$ with respect to the scalar product of $H_{2-r,\overline{v}}$, so

$$(H, \tilde{f})^R = (K_{-r}F, \tilde{f})^R = (F, \Lambda_{2-r}\tilde{f})^R.$$

Since $f = y^{-\frac{2-r}{2}}\tilde{f}$ is in $S_{2-r}(\Gamma, \overline{v})$ and hence holomorphic we have

$$\Lambda_{2-r}\tilde{f} = \Lambda_{2-r}(y^{\frac{2-r}{2}}f) = (z - \overline{z})\frac{\partial f}{\partial \overline{z}} = 0,$$

and therefore $(H, \tilde{f})^R = 0$. □

Theorem 4 now follows from Proposition 6.

*Proof* (of Theorem 4 and of Theorem 2 for $2 - r \neq 1$.)

Let $\phi \in Z^1_{r,v}(\Gamma, \mathcal{P})$ and $g$ and $G$ be constructed as in (23) and (24). In the case $2 - r > 0$ suppose additionally that $(f, \phi) = 0$ for all $f \in S_{2-r}(\Gamma, \overline{v})$. By Lemma 4 in the case $2 - r > 0$ and Remark 4 in the case $2 - r \leq 0$, $G$ satisfies condition (i) of Proposition 6. Hence there is a smooth $F \in H_{-r,\overline{v}}$ and a linear combination of Eisenstein series $E(z) = \sum_{i=1}^{m^*} a_i E_{-r,\overline{v}}^{q_i}(z, \frac{r}{2})$, with

$$G = K_{-r}F + K_{-r}E = K_{-r}(F + E).$$

As stated in Proposition 6, $E$ is only non-zero if $0 \leq 2 - r < 1$, and in this case the Eisenstein series $E_{-r,\overline{v}}^{q_i}(\cdot, \frac{r}{2})$ are smooth functions that grow at most polynomially at each cusp of $\Gamma$. Since $F$ is smooth and in $H_{-r,\overline{v}}$, $F$ also grows at most polynomially at each cusp and so the same is true for $D = E + F$. We have

$$G(z) = y^{\frac{r+2}{2}}\overline{\frac{\partial g}{\partial \overline{z}}(z)} = 2iy\frac{\partial D}{\partial z} - \frac{r}{2}D = 2iy^{\frac{r+2}{2}}\frac{\partial}{\partial z}(y^{-\frac{r}{2}}D).$$

Dividing by $y^{\frac{r+2}{2}}$ and taking the complex conjugate of both sides, we arrive at

$$\frac{\partial g}{\partial \overline{z}}(z) = \frac{\partial}{\partial \overline{z}}(-2iy^{-\frac{r}{2}}\overline{D})(z). \tag{29}$$

Since $D$ is invariant under $|_{-r,\overline{v}}^R$, $\overline{D}$ is invariant under $|_{r,v}^R$. By Lemma 3, the function $\tilde{D}(z) = -2iy^{-\frac{r}{2}}\overline{D}$ is invariant under $|_{r,v}$. This invariance implies that $\tilde{g} = g - \tilde{D}$ satisfies $\tilde{g}|_{r,v}\gamma - \tilde{g} = \phi(\gamma)$ for all $\gamma \in \Gamma$. Since $\tilde{D}$ grows at most polynomially at the cusps of $\Gamma$, $\tilde{g}$ satisfies the growth conditions for functions in $\tilde{\mathcal{Q}}$. Proposition 4 now tells us that $\tilde{g} \in \mathcal{Q}$. Note also that equation (29) implies that $\tilde{g}$ is holomorphic, so $\tilde{g} \in \mathcal{P}$. We finally conclude that $\phi$ is indeed a coboundary in $Z^1_{r,v}(\Gamma, \mathcal{P})$.

The proof above shows in particular that for $2 - r \leq 0$ every cocycle in $Z^1_{r,v}(\Gamma, \mathcal{P})$ is a coboundary and hence $H^1_{r,v}(\Gamma, v) = \{0\}$. This proves Theorem 2 for $2 - r \leq 0$, since $S_{2-r}(\Gamma, \bar{v})$ is also $\{0\}$ in this case. □

*Remark 5* The proof fails if $2 - r = 1$, because Proposition 6 is not available in that case. The only point where we need the assumption $2 - r \neq 1$ in the proof of that proposition is when we show that $\tilde{F}_2$ is in the image of $K_{-r}$, in particular for the construction of the functions $F_2^i \in H_{-r,\bar{v}}$ in (25). The crucial consequence of Proposition 6 is that $G$ is in the image of $K_{-r}$. In the case $2 - r = 1$, we only obtain

$$G = K_{-1}F + \sum_{i=1}^{m^*} \frac{1}{4\pi} \int_{-\infty}^{\infty} \left( E_{1,\bar{v}}^{q_i}\left(\cdot, \frac{1}{2} + i\rho\right), G \right)^R E_{1,\bar{v}}^{q_i}\left(z, \frac{1}{2} + i\rho\right) d\rho.$$

In the notation of the proof of Proposition 6 we have $F = F_1$ and $E = 0$ since $r = 1$. To prove Theorem 2 in this case, one would need to show that the second summand above is in the image of $K_{-1}$.

## 4 Vector-valued modular forms

In this section, we generalise Theorem 2 to vector-valued cusp forms. Let $\rho : \Gamma \to U(n)$ be a unitary representation of $\Gamma$ on $\mathbb{C}^n$ and $v$ a unitary multiplier system of weight $r$. Let $F$ be a function from $\mathcal{H}$ to $\mathbb{C}^n$. The slash operator $|_{\rho,v,r}$ is defined by

$$F|_{r,v,\rho}\gamma(z) = j(\gamma, z)^{-r} \overline{v(\gamma)} \rho(\gamma)^{-1} F(\gamma z).$$

**Definition 10** A function $f : \mathcal{H} \to \mathbb{C}^n$ is a *modular form* for $\Gamma$ of weight $r$, representation $\rho$ and multiplier system $v$ if the following conditions are satisfied:

(i) $f$ is holomorphic on $\mathcal{H}$.
(ii) $f(z) = f|_{r,v,\rho}\gamma(z)$ for all $\gamma \in \Gamma$ and $z \in \mathcal{H}$.
(iii) If $q$ is a cusp of $\Gamma$ and $A\infty = q$, then for any $\epsilon > 0$

$$j(A, z)^{-r} f(Az) \text{ is bounded for } y \geq \epsilon.$$

If $f$ satisfies the additional condition

(iii') If $q$ is a cusp of $\Gamma$ and $A\infty = q$, then there exists an $\epsilon > 0$ such that

$$j(A, z)^{-r} f(Az) = \mathcal{O}_{y \to \infty}(e^{-\epsilon y}),$$

it is a *cusp form*. The set of modular forms or cusp forms of this kind is denoted by $M_r(\Gamma, v, \rho)$ and $S_r(\Gamma, v, \rho)$, respectively.

Let $\mathcal{P}^n$ be the set of vector-valued functions $f(z) = (f_1(z), \ldots, f_n(z))$ such that all $f_i$ are in $\mathcal{P}$. The slash operator $|_{r,v,\rho}$ defines a $\Gamma$-action on $\mathcal{P}^n$ and so we can define the cohomology groups $H^1_{r,v,\rho}(\Gamma, \mathcal{P}^n)$ and $\tilde{H}^1_{r,v,\rho}(\Gamma, \mathcal{P}^n)$. Just as in the one-dimensional case, they turn out to be the same. The proof of this fact relies on a generalisation of Corollary 1:

**Proposition 7** *Let $U \in U(n)$, $s \in \mathbb{R} \setminus \{0\}$ and $g \in \mathcal{P}^n$. Then there exists an $f \in \mathcal{P}^n$ such that*

$$U^* f(z+s) - f(z) = g(z), \quad \forall z \in \mathcal{H}. \tag{30}$$

*Proof* Since $U$ is diagonalisable, there exists a $V \in U(n)$ and a diagonal matrix $D \in U(n)$ with

$$U = V^* D V.$$

Multiplying Eq. (30) by $V$, we get

$$D^* V f(z+s) - V f(z) = V g(z). \tag{31}$$

Let $\epsilon_1, \ldots, \epsilon_n$ be the diagonal entries of $D$ and $G = Vg = (G_1, \ldots, G_n) \in \mathcal{P}^n$. We can use Corollary 1 to find solutions $F_i \in \mathcal{P}$ for

$$\overline{\epsilon_i} F_i(z+s) - F_i(z) = G_i(z).$$

Then $f = V^{-1}(F_1, \ldots, F_n)$ is in $\mathcal{P}^n$ and satisfies (31). □

This can be used to show the following.

**Theorem 8** *Every cocycle in $Z^1_{v,\rho}(\Gamma, \mathcal{P}^n)$ is parabolic.*

### 4.1 Petersson inner product

Let $2 - r > 0$ and $f, g$ be in $S_{2-r}(\Gamma, \overline{v}, \rho^{-1})$. The Petersson inner product of $f$ and $g$ is defined by

$$(f, g) = \int_{\mathcal{F}} \langle f(z), g(z) \rangle y^{-r} dx dy,$$

where $\langle (a_i), (b_i) \rangle = \sum_{i=1}^n a_i \overline{b_i}$ is the usual scalar product on $\mathbb{C}^n$. We will repeat the constructions of Sect. 2.

**Lemma 5** *Let $g$ be in $S_{2-r}(\Gamma, \overline{v}, \rho^{-1})$, then*

$$\phi_g^\infty(z) : \gamma \mapsto \phi_{g,\gamma}^\infty(z) = \overline{\left[ \int_{\gamma^{-1}\infty}^\infty g(\tau)(\tau - \overline{z})^{-r} d\tau \right]}$$

*is a cocycle in $Z^1_{v,\rho}(\Gamma, \mathcal{P}^n)$.*

Again we can use Stokes' theorem to show

$$(f, g) = -C_{2-r} \sum_{m=1}^{n} \int_{A_{i_m}}^{A_{i_m+1}} \left\langle f(z), \overline{\phi_{g,\alpha_{i_m}}^{\infty}(z)} \right\rangle dz.$$

Using this we define a pairing between $S_{2-r}(\Gamma, \overline{v}, \rho^{-1})$ and $H^1_{r,v,\rho}(\Gamma, \mathcal{P}^n)$ as follows. Let $f \in S_{2-r}(\Gamma, \rho^{-1}, \overline{v})$ and $[\phi] \in H^1_{r,v}(\Gamma, \mathcal{P}^n)$ be represented by $\phi$. Then

$$(f, [\phi]) = (f, \phi) = -C_{2-r} \sum_{m=1}^{n} \int_{A_{i_m}}^{A_{i_m+1}} \left\langle f(z), \overline{\phi(\alpha_{i_m})(z)} \right\rangle dz$$

is well defined (independent of the representative $\phi$), and furthermore we have the following theorem, analogous to Theorem 2.

**Theorem 9** *Let $v$ and $\rho$ be as above and $0 < 2 - r \neq 1$. The pairing defined above is perfect, so the map $f \mapsto \phi_f^{\infty}$ induces an isomorphism*

$$S_{2-r}(\Gamma, \overline{v}, \rho^{-1}) \cong H^1_{r,v,\rho}(\Gamma, \mathcal{P}^n).$$

*If $2 - r \leq 0$, we have*

$$S_{2-r}(\Gamma, \overline{v}, \rho^{-1}) \cong H^1_{r,v,\rho}(\Gamma, \mathcal{P}^n) \cong \{0\}.$$

*Proof* All the constructions of Sect. 3 work in the vector-valued case. In particular every statement we cited from [21] is already formulated for vector-valued functions. The fact that every vector-valued modular form of negative weight is 0 is also stated in [21] as a consequence of Satz 5.3, and this generalises the main theorem of [11]. It is also shown that a vector-valued modular form of weight 0 is constant. □

**Acknowledgments** I am thankful to N. Diamantis for suggesting this topic to me and Y. Petridis for a helpful discussion on the proof of Proposition 6. Further thanks are due to F. Strömberg, T. Vavasour and the anonymous referee for a careful reading of this article and many helpful suggestions. I am particularly grateful for the countless comments, corrections and improvements that R. Bruggeman provided during the completion of this article.